\documentclass[reqno,12pt,oneside,a4paper]{amsart} 

\usepackage{bbm,enumerate,amsmath,amstext,mathrsfs,a4}
\usepackage[english]{babel}
\usepackage[latin1]{inputenc}

\newcommand{\R}{\mathbbm{R}}

\newcommand{\Z}{\mathbbm{Z}}

\newcommand{\thh}{\mathrm{Th}}

\newcommand{\F}{\mathbbm{F}}

\newcommand{\IM}{\mathrm{Im}}

\renewcommand{\phi}{\varphi}
\renewcommand{\epsilon}{\varepsilon}
\newcommand{\St}{\mathcal{P}}

\renewcommand{\AA}{\mathscr{A}}

\newcommand{\id}{\mathrm{id}}

\newcommand{\SO}{\mathrm{SO}}
\newcommand{\cl}{\mathrm{cl}}
\newcommand{\Sq}{\mathrm{Sq}}

\usepackage{amsthm}
\theoremstyle{plain}
\newtheorem{thm}{Theorem}[section]
\newtheorem{prop}[thm]{Proposition}
\newtheorem{lem}[thm]{Lemma}

\newtheorem{cor}[thm]{Corollary}
\newtheorem{conj}[thm]{Conjecture}

\theoremstyle{definition}

\newtheorem{defn}[thm]{Definition}

\theoremstyle{remark}

\numberwithin{equation}{section}

\usepackage[all]{xy}
\CompileMatrices

\begin{document}

\title{Secondary characteristic classes of surface bundles}
\author{Søren Galatius}
\address{Aarhus University, Aarhus, Denmark}
\email{galatius@imf.au.dk}

\begin{abstract}
  The Miller-Morita-Mumford classes associate to an oriented surface
  bundle $E\to B$ a class $\kappa_i(E) \in H^{2i}(B;\Z)$.  In this
  note we define for each prime $p$ and each integer $i\geq 1$ a
  secondary characteristic class $\lambda_i(E) \in
  H^{2i(p-1)-2}(B;\Z)/\Z\kappa_{i(p-1)-1}$.  The mod $p$ reduction
  $\lambda_i(E) \in H^*(B; \F_p)$ has zero indeterminacy and satisfies
  $p\lambda_i(E) = \kappa_{i(p-1)-1}(E) \in H^*(B;\Z/p^2)$.
\end{abstract}

\maketitle

\section{Introduction and statement of results}

Recall that any bundle $\pi:E\to B$ of oriented surfaces with finite
dimensional base $B$ has an embedding $j: E \to B\times \R^{N+2}$ over
$B$.  For $N$ large, $j$ is unique up to isotopy.  A choice of
embedding $j$ induces a transfer map
\begin{equation*}
  \xymatrix{
    {B_+ \wedge S^{N+2}} \ar[r]^-{\pi_!} & {\thh(\nu j)}
  }
\end{equation*}
The embedding $j: E \to B\times \R^{N+2}$ also induces classifying
maps
\begin{equation*}
  \xymatrix{
    {T_\pi E \ar[r]^-{\cl(T_\pi E)}} \ar[d] & {\SO(N+2)
      \times_{\SO(N)\times \SO(2)} \R^2} \ar[d] \\
    E \ar[r] & {\SO(N+2) / \SO(N)\times \SO(2)}\\
  }    
\end{equation*}
and
\begin{equation*}
  \xymatrix{
    {\nu j \ar[r]^-{\cl(\nu j)}} \ar[d] & {\SO(N+2)
      \times_{\SO(N)\times \SO(2)} \R^N} \ar[d] \\
    E \ar[r] & {\SO(N+2) / \SO(N)\times \SO(2)}\\
  }    
\end{equation*}

For brevity, write $U = U_N = \SO(N+2)\times_{\SO(N)\times \SO(2)}
\R^2$ and $U^\perp = U_N^\perp = \SO(N+2)\times_{\SO(N)\times \SO(2)}
\R^N$.  We get the composition
\begin{equation*}
  \alpha = \thh(\cl(\nu j)) \circ \pi_!: B_+ \wedge S^{N+2} \to
  \thh(U_N^\perp)
\end{equation*}
Recall that there is a Thom class $\lambda_{U^\perp} \in
H^N(\thh(U^\perp),*;\Z)$ and that we have $H^{N+*}(\thh(U^\perp),*;\Z)
= \Z[e(U)] .  \lambda_{U^\perp}$ for $* < N$.  The definition of the
$\kappa$-classes is
\begin{equation*}
  \kappa_i E = \alpha^*(e(U)^{i+1} . \lambda_{U^\perp}) = \pi_!^*(
  e(T_\pi E)^{i+1} . \lambda_{\nu j}) \in H^{2i}(B;\Z)
\end{equation*}

In this paper we define secondary characteristic classes of surface
bundles.  The definition involves Toda brackets.  In
section~\ref{sec:second-comp} we recall some generalities about Toda
brackets.  By a surface bundle we shall mean a fibre bundle with
closed oriented smooth two-dimensional fibres.
\begin{lem}\label{lem:0}
  Let $p$ be a prime, and let $\St^i$ denote the Steenrod power
  operation.  When $p=2$, write $\St^i = \Sq^{2i}$ and $\beta\St^i =
  \Sq^{2i+1}$.  Given a surface bundle $\pi: E\to B$, let $\alpha: B_+
  \wedge S^{N+2} \to \thh(U_N^\perp)$ be as before and let $\lambda:
  \thh(U_N^\perp) \to K(\Z,N)$ be the Thom class.  Then the Toda
  bracket
  \begin{equation*}
    \{\beta\St^i, \lambda, \alpha\} \subseteq
    H^{2i(p-1)-2+N}(B_+\wedge S^{N+2};\Z) = H^{2i(p-1)-2}(B;\Z)
  \end{equation*}
  is defined with indeterminacy $\Z\kappa_{i(p-1)-1}$.
\end{lem}
\begin{defn}
  With notation as in Lemma~\ref{lem:0} define
  \begin{equation*}
    \lambda_i(E) = (-1)^i \{\beta \St^i, \lambda, \alpha\} \in
    H^{2i(p-1)-2} (B;\Z) / \Z\kappa_{i(p-1)-1}
  \end{equation*}
\end{defn}
\begin{thm}
  \label{thm:A}
  The mod $p$ reduction $\lambda_i(E) \in H^*(B;\F_p)$ has
  zero indeterminacy and satisfies
  \begin{equation*}
    p\lambda_i(E) = \kappa_{i(p-1)-1} \in H^*(B;\Z/p^2)
  \end{equation*}
  More generally we have the following in integral cohomology
  \begin{equation*}
    \kappa_{i(p-1)-1} \in p\lambda_i(E)
  \end{equation*}
\end{thm}
\begin{thm}\label{thm:B}
  \begin{enumerate}[(i)]
  \item If $\pi:E\to B$ and $\pi': E' \to B$ are surface bundles, then
    \begin{equation*}
      \lambda_i(E\amalg E') = \lambda_i(E) + \lambda_i(E')
    \end{equation*}
  \item If $\pi: E\to B$ is a surface bundles and $\pi': E'\to B$ is
    obtained from $E$ by fibrewise surgery, then
    \begin{equation*}
      \lambda_i E = \lambda_i E'
    \end{equation*}
  \item If $\pi: E \to B$ and $\pi': E' \to B$ are bundles of compact,
  non-closed surfaces with $\partial E = S^1 \times B = \partial E'$,
  then
  \begin{equation*}
    \lambda_i(E \cup_{S^1\times B} E') = \lambda_i(E \cup_{S^1\times
    B} (D^2 \times B)) + \lambda_i(E' \cup_{S^1\times B} (D^2 \times B))
  \end{equation*}
  \end{enumerate}
\end{thm}

As an application of secondary classes we prove the following
strengthening of a theorem of~\cite{GMT}:
\begin{thm}\label{thm:C}
  Let $p$ be a prime and $s\geq 1$.  Then the reduction of
  $\kappa_{ps(p-1)-1}$ mod $p^2$ vanishes:
  \begin{equation*}
    \kappa_{ps(p-1)-1} = 0 \in H^*(B; \Z/p^2)
  \end{equation*}
\end{thm}

Theorem~\ref{thm:C} proves part of the following conjecture.
\begin{conj}\label{conj:D}
  Let $s \geq 1$ and $v\geq 0$.  Then
  \begin{equation*}
    \kappa_{p^vs(p-1)-1} = 0 \in H^*(B;\Z/p^{v+1})
  \end{equation*}
\end{conj}
If the conjecture is true, then $\kappa_{p^vs(p-1)-1}$ can be divided
by $p^{v+1}$.  In~\cite{GMT} we prove that this holds modulo torsion.
It is also proved in \cite{GMT} that the statement of
Conjecture~\ref{conj:D} is best possible in the sense that if
$s\not\equiv 0\pmod{p}$, then $\kappa_{p^vs(p-1)-1} \neq 0 \in
H^*(B;\Z/p^{v+2})$.  I hope to return to
Conjecture~\ref{conj:D} at a later time.

\section{Secondary composition}\label{sec:second-comp}
We recall the definition of secondary compositions (Toda brackets).
For further details see \cite{Toda}.

All spaces and maps are pointed.  The reduced suspension $SX$ is
regarded as the pushout of $\xymatrix{X \wedge [-1,0] & X \ar[l]
  \ar[r] & X \wedge [0,1]}$ where $-1 \in [-1,0]$ and $1\in [0,1]$ are
the basepoints.  Thus, two nullhomotopies $F: X\wedge [-1,0] \to Y$
and $G: X\wedge[0,1] \to Y$ induce a map $G-F: SX \to Y$.

For a sequence of maps
\begin{equation*}
  \xymatrix{
    X \ar[r]^f & Y\ar[r]^g & Z\ar[r]^h & W
  }
\end{equation*}
with $g\circ f \simeq 0$ and $h \circ g \simeq 0$, a choice of
null-homotopies $F: g\circ f \simeq 0$ and $G: h\circ g \simeq 0$
determines a map
\begin{equation*}
  h\circ F - G \circ (f\wedge [-1,0]) : SX \to W
\end{equation*}
We define the \emph{secondary composition} to be the subset $\{h,g,f\}
\subseteq [SX,W]$ of homotopy classes of maps obtained in this
fashion, as $F,G$ ranges over all nullhomotopies.

Recall that $[SX,W] = [X,\Omega W]$ is a group.
\begin{lem}
  \label{lem:1}
  $\{h,g,f\}$ depends only on the homotopy classes of $h$, $g$, and
  $f$.  If $\{h,g,f\}$ is defined, then it gives a unique element in
  the double coset,
  \begin{equation*}
    \{h,g,f\} \in h\circ[SX,Z]\setminus[SX,W]/[SY,W]\circ Sf
  \end{equation*}
  If $[SX,W]$ is abelian, then
  \begin{equation*}
    \{h,g,f\} \in [SX,W]/\big( h\circ[SX,Z] + [SY,W]\circ Sf\big)
  \end{equation*}
\end{lem}
\begin{proof}
  See \cite[Lemma 1.1]{Toda}.
\end{proof}
\begin{prop}\label{prop:2}
  For a sequence of maps
  \begin{equation*}
    \xymatrix{
      X \ar[r]^f & Y\ar[r]^g & Z\ar[r]^h & W \ar[r]^k & V
    }
  \end{equation*}
  we have
  \begin{enumerate}[(i)]
  \item $\{k,h,g\} \circ f \subseteq \{k,h,g\circ f\}$
  \item $\{k,h,g\circ f\} \subseteq \{k,h\circ g, f\}$
  \item $\{k\circ h, g, f\} \subseteq \{k,h\circ g, f\}$
  \item $k\circ\{h,g,f\} \subseteq \{k\circ h, g,f\}$
  \end{enumerate}
\end{prop}
\begin{proof}
  See \cite[Proposition 1.2]{Toda}.
\end{proof}
\begin{prop}\label{prop:divp}
  Let 
  \begin{equation*}
    \xymatrix{
      {K(\Z,n)}\ar[r]^p &
      {K(\Z,n)}\ar[r]^{\rho} &
      {K(\F_p,n)}\ar[r]^{\beta} &
      {K(\Z,n+1)}
    }
  \end{equation*}
  represent multiplication by $p$, reduction mod $p$, and the mod $p$
  Bockstein, respectively.  Then
  \begin{equation*}
    \id \in \{\beta, \rho, p\} \subseteq [SK(\Z,n), K(\Z,n+1)] =
    [K(\Z,n), K(\Z,n)]
  \end{equation*}
  \qed
\end{prop}
\begin{cor}\label{cor:divp}
  Let $c: X\to K(\Z,n)$ represent a cohomology class.  Let $\rho$ and
  $\beta$ be as in Proposition~\ref{prop:divp}.  Then
  \begin{equation*}
    \{\beta, \rho, c\} = \tfrac1p c + \Z c \subseteq H^n(X) =
    [SX,K(\Z,n+1)]
  \end{equation*}
  where
  \begin{equation*}
    \tfrac1p c = \{c' \vert pc' = c\}
  \end{equation*}
\end{cor}
\begin{proof}
  Clearly the two sides have the same indeterminacy $\Z c + \beta
  H^{n-1}(X;\F_p)$, so all we need to check is that if $pc' = c$, then
  $c'\in \{\beta, \rho, c\}$.  But this follows from
  Proposition~\ref{prop:divp}:
  \begin{equation*}
    \{\beta, \rho, p\circ c'\} \supseteq \{\beta, \rho, p\} \circ c'
    \ni c'
  \end{equation*}
\end{proof}

\section{Elementary properties of the secondary
  classes}\label{sec:defin-second-class} 

Consider the oriented Grassmannian $\SO(N+2)/\SO(N)\times\SO(2)$.  Let
$U = U_N = \SO(N+2)\times_{\SO(N)\times\SO(2)} \R^2$ be the canonical
oriented 2-dimensional vectorbundle and let $U^\perp = U_N^\perp =
\SO(N+2)\times_{\SO(N)\times\SO(2)} \R^N$ be its orthogonal
complement.
\begin{lem}[\cite{GMT}]
  \label{lem:GMT}
  In $H^*(\thh(U^\perp),* ; \F_p)$ we have that
  \begin{equation*}
    \St^i \lambda_{U^\perp} = (-1)^i e^{i(p-1)} \lambda_{U^\perp}
  \end{equation*}
\end{lem}
\begin{proof}
  Let $\St = \sum_i \St^i$.  Then $\St(\lambda_U) =
  (1+e(U)^{p-1})\lambda_U$.  Since $\lambda_{U\oplus U^\perp} =
  \lambda_U \cup \lambda_{U^\perp}$ we get
  \begin{equation*}
    \lambda_U\cup \lambda_{U^\perp} = \lambda_{U\oplus U^\perp} =
    \St(\lambda_{U\oplus U^\perp}) = \St(\lambda_U) \cup
    \St(\lambda_{U^\perp}) = (1+e(U)^{p-1})\lambda_U\cup
    \St(\lambda_{U^\perp})
  \end{equation*}
  and hence
  \begin{equation*}
    \St(\lambda_{U^\perp}) = (1+e(U)^{p-1})^{-1} \lambda_{U^\perp} =
    \bigg( \sum_i (-1)^i e(U)^{i(p-1)} \bigg) \lambda_{U^\perp}
  \end{equation*}
\end{proof}
\begin{proof}[Proof of Lemma~\ref{lem:0}]
  Clearly $l\circ \alpha\simeq 0$.  The cohomology of the Grassmannian
  $\SO(N+2)/\SO(N)\times\SO(2)$ vanishes in odd degrees (when $N$ is
  larger than the degree), so $\beta\St^i \circ \lambda \simeq 0$.
  Therefore $\{\beta\St^i, \lambda, \alpha\}$ is defined.  It follows
  from Lemma~\ref{lem:1} that the indeterminacy is
  $\Z\kappa_{i(p-1)-1}$.
\end{proof}

\begin{proof}[Proof of Theorem~\ref{thm:A}]
  This follows from Proposition \ref{prop:2} and
  Corollary~\ref{cor:divp} and the diagram:
  \begin{equation*}
    \xymatrix{
      {B_+ \wedge S^{N+2}} \ar[r]^-\alpha \ar[rd]_-{\kappa_{i(p-1)-1}} &
      {\thh(U_N^\perp)} \ar[r]^-\lambda \ar[d]^-{e^{i(p-1)}\lambda} &
      {K(\Z, N)}\ar[d]^-{\St^i}\\
      & {K(\Z,N+2i(p-1))} \ar[r]^-\rho & 
      {K(\F_p, N+2i(p-1))} \ar[d]^\beta\\
      && {K(\Z,N+2i(p-1)+1)}\\
    }
  \end{equation*}
  Indeed, Proposition \ref{prop:2} gives the inclusions
  \begin{align*}
    \{\beta, \rho, \kappa_{i(p-1)-1}\} &= \{\beta, \rho,
    (e^{i(p-1)}\lambda)\circ \alpha\} \subseteq \{\beta, \rho\circ
    (e^{i(p-1)}\lambda), \alpha\}\\
    &= (-1)^i\{\beta, \St^i \lambda, \alpha\} \supseteq (-1)^i\{\beta\St^i,
    \lambda, \alpha\} = \lambda_i(E).
  \end{align*}
  Then Lemma \ref{lem:1} proves that the first inclusion is an
  equality since the two sides have the same indeterminacy $\IM(\beta)
  + \Z\kappa_{i(p-1)-1}$.  Therefore by Corollary~\ref{cor:divp}
  \begin{equation*}
    \lambda_i(E) \subseteq \{\beta, \rho, \kappa_{i(p-1)-1}\} =
    \tfrac1p\kappa_{i(p-1)-1} + \Z \kappa_{i(p-1)-1},
  \end{equation*}
  and hence
  \begin{equation*}
    p\lambda_i(E) \subseteq (1+p\Z) \kappa_{i(p-1)-1}.
  \end{equation*}
  Since they have the same indeterminacy, they are equal.
\end{proof}
\begin{proof}[Proof of Theorem~\ref{thm:B}]
  \emph{(i)} follows from the additivity of $\alpha$, i.e. the
  property that $\alpha(E\amalg E') = \alpha(E) + \alpha(E') \in
  [B_+\wedge S^{N+2},\thh(U_N^\perp)]$.  Similarly \emph{(ii)} follows
  from the property that $\alpha(E) = \alpha(E')$ when $E'$ is
  obtained from $E$ by fibrewise surgery.  \emph{(iii)} follows from
  \emph{(i)} and \emph{(ii)} since $E\cup_{S^1\times B} E'$ is
  obtained from $\big(E\cup_{S^1\times B} (D^2\times B)\big) \amalg
  \big(E\cup_{S^1\times B} (D^2\times B)\big)$ by fibrewise surgery.
\end{proof}
\section{A variant of $\lambda_{ps}$}
The goal of this section is to prove Theorem~\ref{thm:C}.  The
definition and properties of $\lambda_i$ proves that $\kappa_{i(p-1)}$
is divisible by $p$.  When $i = ps$, a variant of $\lambda_{ps}$ can
be used to prove that $\kappa_{ps(p-1)-1}$ is divisible by $p^2$.
\begin{defn}
  Let $s \geq 0$ and consider the Steenrod algebra $\AA_p$.  When
  $p=2$ we write $\St^i = \Sq^{2i}$ and $\beta\St^i = \Sq^{2i+1}$ as
  before.  Define $\theta_s\in \AA_p$ by
  \begin{equation*}
    \theta_s = \sum_{j=0}^s (-1)^j \binom{(p-1)(s-j)}j
    \St^{ps-j}\St^{j} = \St^{ps} + \text{terms of lenght 2}
  \end{equation*}
  Define vectors $v_s, w_s\in\AA_p$ by
  \begin{equation*}
    w_s = (\St^0, \dots \St^s), \quad v_s = (\St^{ps}, \dots,
    (-1)^j\binom{(p-1)(s-j)-1}j \St^{ps-j}, \dots, \St^{(p-1)s}).
  \end{equation*}
\end{defn}
\begin{lem}
  \begin{enumerate}[(i)]
  \item In $H^*(\thh(U^\perp),*;\F_p)$ we have that $\theta_s\lambda_{U^\perp} =
  e^{ps(p-1)}\lambda_{U^\perp}$.
  \item $v_s^T\beta w_s = \beta\theta_s$.
  \end{enumerate}
\end{lem}
\begin{proof}
  \emph{(i)} This is similar to Lemma~\ref{lem:GMT}, using the fact
  that the admissible terms of length 2 act trivially on
  $\lambda_{U^\perp}$. Formula \emph{(ii)} is the Adem relation for
  $\St^{(p-1)s}\beta\St^s$.
\end{proof}
\begin{defn}
  Let $\alpha, \lambda, \theta_s$ be as above.  Define the secondary
  characteristic class
  \begin{equation*}
    \tilde \lambda_{ps} = (-1)^s\{\beta\theta_s, \lambda, \alpha\} \in
    H^{2ps(p-1)-2}(B, \Z)/\Z\kappa_{ps(p-1)-1}
  \end{equation*}
\end{defn}
Notice that $\tilde\lambda_{ps}$ satisfies the same formal properties
as $\lambda_{ps}$.  In particular $p\tilde\lambda_{ps} =
(1+p\Z)\kappa_{ps(p-1)-1}$.  In general $\tilde \lambda_{ps} \neq
\lambda_{ps}$.
\begin{proof}[Proof of Theorem~\ref{thm:C}]
  We have
  \begin{align*}
    (-1)^s\rho\circ\{\beta\theta_s, \lambda, \alpha\} \subseteq
    (-1)^s\{\rho\circ\beta\theta_s, \lambda, \alpha\} & = 
    (-1)^s\{v_s^T\beta w_s, \lambda, \alpha\} \\ &\supseteq (-1)^sv_s^T\{\beta w_s,
    \lambda, \alpha\} 
  \end{align*}
  and it is seen that all the inclusions are equalities since the
  indeterminacy vanishes.  Since
  \begin{equation*}
    (-1)^s\{\beta w_s, \lambda, \alpha\} \in \prod_{i=0}^s H^{N+2i(p-1)}(B_+
    \wedge S^{N+2}; \F_p) = \prod_{i=0}^s H^{2i(p-1)-2}(B; \F_p),
  \end{equation*}
  $v^T$ will vanish since $H^*(B;\F_p)$ is an unstable $\AA_p$-module.

  Hence the mod $p$ reduction of $\tilde\lambda_{ps}$ vanishes, so
  $\kappa_{ps(p-1)-1} = p\tilde\lambda_{ps} = 0 \in H^*(B;\Z/p^2)$.
\end{proof}

\end{document}